\documentclass[12pt]{article}
\usepackage{fancyhdr}
\usepackage{graphics,epstopdf}
\usepackage{float}
\usepackage{courier}
\usepackage{epsfig}
\usepackage{lipsum}
\usepackage{amssymb,amsmath}
\usepackage{amsfonts}
\usepackage{color}
\usepackage{todonotes,dsfont}
\RequirePackage{lineno}

\newtheorem{theorem}{Theorem}[section]

\newtheorem{proposition}[theorem]{Proposition}

\newtheorem{remark}[theorem]{Remark}
\newtheorem{definition}[theorem]{Definition}
\newtheorem{example}[theorem]{Example}

\catcode`@=11 \@addtoreset{equation}{section}
\renewcommand\theequation{\thesection.\@arabic\c@equation}
\catcode`@=12

\makeatother

\begin{document}

\title{\bf Multiframelet Properties on $\mathbb{Q}_{p}$}

\author{Debasis Haldar\footnote{E-mail : debamath.haldar@gmail.com} \\ \\
Department of Mathematics, National Institute of Technology Rourkela,\\ Rourkela - 769008, India}

\date{}
\maketitle

\begin{abstract}
  This paper produces various results on $p$-adic multiframelet. Multiframelet is a frame-like sequence generated by multiple functions along with wavelet structure. Various properties of multiframelet in $L^{2}(\mathbb{Q}_{p})$ have been analyzed. Also multiframelet operator on $p$-adic setting has been produced and characterized. Furthermore, multiframelet set in $\mathbb{Q}_{p}$ has been engendered and scrutinized.
\end{abstract}



\section{Introduction}

Frame is a sequence which allows to write every element of the Hilbert space as linear combination of corresponding frame elements. Coefficients of this representation are not unique and that is why frames are sometimes said to be overcomplete system. Its primary objective was substitution or extension of Riesz bases or orthonormal bases in Hilbert spaces. Duffin and Schaeffer were introduced frame in \cite{1} for computing the coefficients in a linear combination of the vectors of a linearly dependent spanning set. After several decades, Daubechies, Grossmann and Meyer \cite{3} took the key steps of connecting frames with wavelets and Gabor systems. Heil \cite{11} devoloped frame theory for general Hilbert space. Grochenig \cite{8} gave the nontrivial extension of frames to Banach spaces. Frames, having wavelet structures, have been popularized through several generalizations and significant applications in signal processing, image processing etc, for detail discussion regarding the same we refer \cite{Bh18,9,DDRR19,10,3,1,HS19,12,RDD17,2}. \\

This paper is arranged as follows; In Section 2, we have discussed about basics of $p$-adic numbers. Multiframelet and its various properties on $p$-adic setting have discussed in Section 3. Section 4 has devoted for multiframelet operator and its dual and also representation of arbitrary elements using multiframelet operator. Furthermore, multiframelet set has carried out in Section 5.

\section{Preliminaries}

The field $\mathbb{Q}_{p}$ of $p$-adic numbers is defined as the completion of $\mathbb{Q}$ with respect to metric topology induced by the $p$-adic norm $|\cdot|_{p}$. The $p$-adic norm is defined as follows
$$|x|_{p}= \begin{cases}
0 & \text{ if } x=0, \\
p^{- \gamma} & \text{ if } x=p^{\gamma}\frac{m}{n} \neq 0 ,~~ p \not| ~ mn.
\end{cases}$$

This norm has the ultrametric property $\lvert x+y \rvert_{p} \leq \text{max} \{\lvert x \rvert_{p}, \lvert y \rvert _{p} \}$. Here the equality holds if and only if $|x|_{p} \neq |y|_{p}$. Thus the $p$-adic norm is non-Archimedean. Every non-zero $p$-adic number $x$ has canonical representation
$x= \sum \limits_{j=\gamma}^{\infty} x_{j}p^j,$
where $x_{j} \in \{ 0,1,...., p-1 \}$ with $x_{\gamma} \neq 0$ and $\gamma \in \mathbb{Z}$. The fractional part of $x$ is $\{x\}_{p}:=\sum \limits_{j=\gamma}^{-1}x_{j}p^j$. Thus  $\{x\}_{p} =0$ if and only if $\gamma \geq 0$. Further $\{0\}_{p}:=0$. The ring of $p$-adic integers is denoted as $\mathbb{Z}_{p}$ and the set of fractional numbers $I_{p}$ are given by $\mathbb{Z}_{p} = \{ x \in \mathbb{Q}_{p} : \{x\}_{p}=0 \}$ and $I_{p} = \{ x \in \mathbb{Q}_{p}: \{x\}_{p} = x \}$. Therefore, $I_{p}=\{ \frac{a_{-\gamma}}{p^{\gamma}} + \frac{a_{-\gamma +1}}{p^{\gamma -1}} + \ldots + \frac{a_{-1}}{p} \in \mathbb{Q}_{p} : 0 \leq a_{i} \leq p-1, i \in \mathbb{Z}, a_{-\gamma} \neq 0 \}$. The additive character $\chi_{p}$ on the field $\mathbb{Q}_{p}$ is defined by $\chi_{p}(x)= e^{2\pi i \{x\}_{p}},~~~\forall~ x \in \mathbb{Q}_{p}$.\\

A ball $B_{\gamma}(a)$ centered at $a\in \mathbb{Q}_{p}$ and radius $p^\gamma$, where $\gamma \in \mathbb{Z}$, is defined by $B_{\gamma}(a):= \{x \in \mathbb{Q}_{p}: \lvert x-a \rvert_{p} \leq p^{\gamma} \}.$ The ball $B_{\gamma}(a)$ is compact, open and every point of the ball is a center of the ball. Moreover, any two balls are either disjoint or one is contained in the other. $\mathbb{Q}_{p}$ is totally disconnected, locally compact and has no isolated points. There exists Haar measure $dx$ on $\mathbb{Q}_p$ which is positive and translations invariant, i.e., $d(x+a)=dx$, for all $a \in I_{p}$. It is normalized by the equality $\int \limits_{\mathbb{Z}_{p} }dx =1$. Furthermore,
$d(ax)= |a|_{p}~ dx,~~~  \forall~a \in \mathbb{Q}_{p} \setminus \{0 \}.$

The space of all $\mathbb{C}$-valued functions on $\mathbb{Q}_{p}$, square integrable with respect to the measure $dx$, is denoted by $ L^{2}(\mathbb{Q}_{p})$. The inner product in this space is defined by
$$\langle f,g \rangle = \int \limits_{\mathbb{Q}_{p}} f(x) \overline{g(x)}dx,~~~~ \forall~ f, g \in L^{2}(\mathbb{Q}_{p}).$$

The Fourier transform of $f \in L^{2}(\mathbb{Q}_{p})$ is defined as
\begin{eqnarray*}
	\hat{f}(\xi)= \int \limits_{\mathbb{Q}_p}\chi_{p}(\xi x)f(x) dx, ~~~~\forall~ \xi \in \mathbb{Q}_{p}.
\end{eqnarray*}

This Fourier transform satisfy following equality, namely Parseval theorem \cite{T75},
\begin{eqnarray*}
	\langle f,g \rangle = \langle \hat{f}, \hat{g} \rangle ~~~~ \forall ~ f, g \in L^{2}(\mathbb{Q}_{p}).
\end{eqnarray*}

Throughout this paper we denote $\mathcal I$ as the identity operator, $\mathfrak{L}=\{1,2, \cdots, L \}$, $\mathds{1}_{E}$ as characteristic function of the set $E$, $p$ is a prime number, $\mathcal{R}_{\mathcal{U}}$ and $\mathcal{N}_{\mathcal{U}}$ are denoted as range and kernel respectively of bounded linear operator $\mathcal{U}$.

\section{Multiframelet on $\mathbb{Q}_{p}$}

\begin{definition} (Multiframelet).
	A set of functions $\mathfrak{f}=\{f^{(1)}, \ldots, f^{(L)}\} \subset L^{2}(\mathbb{Q}_{p})$ is said to be a multiframelet of order L if $\{f^{(l)}_{j,a}:=p^{\frac{j}{2}}f^{(l)}(p^{-j} \cdot -a) : j \in \mathbb{Z},~ a \in I_{p},~ l \in \mathfrak{L}\}$ is a frame for $L^{2}(\mathbb{Q}_{p})$ i.e. $\exists ~ A,B > 0$ such that
	\begin{equation} \label{fie}
	A\left \| g \right|\|^{2} \leq \sum \limits_{l \in \mathfrak{L}} \sum \limits_{j \in \mathbb{Z}} \sum \limits_{a \in I_{p}}| \langle g,f^{(l)}_{j,a} \rangle |^{2} \leq B\left \| g \right \|^{2}, ~~~~~~~~ \forall ~ g \in L^{2}(\mathbb{Q}_{p}).
	\end{equation}
\end{definition}
When $L=1$, $\mathfrak{f}$ is simply said to be a framelet. $A$ and $B$ are said to be lower and upper bounds respectively for the multiframelet. Clearly, they are not unique. The optimal lower multiframelet bound is the supremum of all lower multiframelet bounds and the optimal upper multiframelet bound is the infimum of all upper multiframelet bounds. A multiframelet which ceases to be a multiframelet on the removal of any one of its vectors is termed an exact multiframelet. A multiframelet which is not exact is called inexact multiframelet. $\mathfrak{f}$ is said to be a tight multiframelet if it is possible to choose $A=B$ and $\mathfrak{f}$ is said to be a normalized tight multiframelet or Parseval multiframelet if it is possible to choose $A=B=1$. In this case, the multiframelet is also said to be $A$-tight multiframelet. Every orthonormal basis is a Parseval multiframelet but a Parseval multiframelet need not be orthogonal or a basis.

\begin{example} \label{koz}
	Kozyrev's multiwavelet \cite{5} of order $(p-1)$ is given by
	\begin{equation*}
	\theta_{k}(x) = \chi_{p}(p^{-1}kx) \mathds{1}_{\mathbb{Z}_p}(x) ,  \hspace{1 cm}  x \in \mathbb{Q}_{p},
	\end{equation*}
	where $ k = 1, 2,\ldots, p-1$ .
\end{example}

\begin{example} \label{ks}
	Khrennikov and Shelkovich's multiwavelet \cite{4} of order $(p-1)p^{m-1}$ is given by
	\begin{equation*}
	\theta ^{(m)}_{s}(x)=\chi_{p}(sx) \mathds{1}_{\mathbb{Z}_p}(x) , \hspace{.3 cm} x \in \mathbb{Q}_{p},
	\end{equation*}
	where $s \in J_{p,m}:= \{ \frac{s_{-m}}{p^m} + \ldots +\frac{s_{-1}}{p} :~ s_{-j} = 0, 1,\ldots, p-1$; $j = 1, 2,\ldots, m;~ s_{-m} \neq 0 \}$, $m \in \mathbb{N}$ is fixed. \end{example}

\begin{remark}
	Example \ref{koz} and Example \ref{ks} are normalized tight multiframelet. Every multiwavelet is also normalized tight multiframelet. Moreover, any orthonormal basis is an exact frame as deleting any term from an orthonormal basis gives 0 for the middle portion in the Equation \eqref{fie} for that particular deleted $g$. This imply $g=0$, which is contradiction to orthonormal basis.
\end{remark}

\noindent From Equation \eqref{fie}, in a tight multiframelet i.e. when $A = B$, we have
\begin{equation*}
\sum \limits_{l \in \mathfrak{L}} \sum \limits_{j \in \mathbb{Z}} \sum \limits_{a \in I_{p}}| \langle g,f^{(l)}_{j,a} \rangle |^{2}=A \left \| g \right \|^{2}, ~~~~~~~~ \forall ~ g \in L^{2}(\mathbb{Q}_{p}).
\end{equation*}
By pulling $\frac{1}{A}$ into the sum, this is equivalent to
\begin{equation}
\sum \limits_{l \in \mathfrak{L}} \sum \limits_{j \in \mathbb{Z}} \sum \limits_{a \in I_{p}}| \langle g,\frac{1}{\sqrt{A}}f^{(l)}_{j,a} \rangle |^{2}= \left \| g \right \|^{2}, ~~~~~~~~ \forall ~ g \in L^{2}(\mathbb{Q}_{p}). \nonumber
\end{equation}
Thus the family $\{ \frac{1}{\sqrt{A}} f_{j,a}^{(l)} : j \in \mathbb{Z},~ a \in I_{p},~ l \in \mathfrak{L} \}$ is a 1-tight frame. So any tight frame can be changes to normalized tight frame. In other word, we can use definition of tight frame as follows
\begin{equation} \label{fe}
\sum \limits_{l\in \mathfrak{L}} \sum \limits_{j \in \mathbb{Z}} \sum \limits_{a \in I_{p}}| \langle g, f^{(l)}_{j,a} \rangle |^{2}= \left \| g \right \|^{2} ~~~~~~~~ \forall ~ g \in L^{2}(\mathbb{Q}_{p}).
\end{equation}
assuming that $\{ f_{j,a}^{(l)} : j \in \mathbb{Z},~ a \in I_{p},~ l\in \mathfrak{L} \}$ has been properly normalized. It follows from the Equation \eqref{fe} that normalized tight frame $\{ f_{j,a}^{(l)} : j \in \mathbb{Z},~ a \in I_{p},~ l\in \mathfrak{L} \}$ is an orthonormal basis of $L^{2}(\mathbb{Q}_{p})$ if and only if $\|f^{(1)}\|=\|f^{(2)}\|= \ldots =\|f^{(L)}\|=1.$

The operator $\mathcal{T} : \ell^{2}(\mathfrak{L} \times \mathbb{Z} \times I_{p}) \rightarrow L^{2}(\mathbb{Q}_{p})$ defined by $c(l,j,a) \mapsto \sum \limits_{l\in \mathfrak{L}} \sum \limits_{j \in \mathbb{Z}} \sum \limits_{a \in I_{p}} c(l,j,a) f^{(l)}_{j,a}$ is called synthesis operator or pre-multiframelet operator. The adjoint operator of $\mathcal{T}$ is $\mathcal{T}^{*} :  L^{2}(\mathbb{Q}_{p}) \mapsto \ell^{2}(\mathfrak{L} \times \mathbb{Z} \times I_{p})$ defined as $g \mapsto \{ \langle g, f_{j,a}^{(l)} \rangle :~ j \in \mathbb{Z},~ a \in I_{p},~ l \in \mathfrak{L}\}$. This is called analysis operator. We now prove some basic results on multiframelet.

\begin{proposition}
	Multiframelets are bounded above with bound $\sqrt{B}$.
\end{proposition}

\noindent \textbf{Proof :}
	For a fixed $l^{\prime}, a^{\prime}, j^{\prime}$; we have
	\begin{eqnarray*}
		&&\|f_{j^{\prime},a^{\prime}}^{(l^{\prime})}\|^{4} = |\langle f_{j^{\prime},a^{\prime}}^{(l^{\prime})}~,~ f_{j^{\prime},a^{\prime}}^{(l^{\prime})}\rangle|^{2} \leq \sum \limits_{l\in \mathfrak{L}} \sum \limits_{j \in \mathbb{Z}} \sum \limits_{a \in I_{p}} |\langle f_{j^{\prime},a^{\prime}}^{(l^{\prime})}~,~ f_{j,a}^{(l)}\rangle|^{2} \leq B \|f_{j^{\prime},a^{\prime}}^{(l^{\prime})}\|^{2} \\ 
		&\text{i.e.}&  \|f_{j^{\prime},a^{\prime}}^{(l^{\prime})}\|^{2} \leq B ~~ \Rightarrow ~~ \|f_{j^{\prime},a^{\prime}}^{(l^{\prime})}\| \leq \sqrt{B}.
	\end{eqnarray*}
	So multiframelets are bounded above.

\begin{proposition}
	Inexact multiframelets are not bases.
\end{proposition}

\noindent \textbf{Proof :}
	Let $\{ f^{(1)}, \ldots, f^{(L)}\}$ be an inexact multiframelet. Then there exists either $j^{\prime} \in \mathbb{Z}$ or $a^{\prime} \in I_{p} $ or $l^{\prime} \in \mathfrak{L}$ such that $\{f_{j,a}^{(l)}\}_{l \neq l^{\prime} ~or~ j \neq j^{\prime} ~or~ a \neq a^{\prime}}$ is a frame and therefore complete. But a basis is a minimal spanning set and hence no subset of a basis is complete. Therefore inexact multiframelets can not be basis.

In the following proposition we present a result on completeness of multiframelet.

\begin{proposition} \label{lower}
	The lower multiframelet condition ensures that a multiframelet is complete i.e. its closed linear span is the whole space $L^{2}(\mathbb{Q}_{p})$.
\end{proposition}

\noindent \textbf{Proof :}
	We will show it by contradiction. Let $V= \overline{span} \{f^{(l)}_{j,a} : j \in \mathbb{Z},~ a \in I_{p},~l \in \mathfrak{L}\}$. Then $L^{2}(\mathbb{Q}_{p})= V \oplus V^{\bot}$. If $V^{\bot} \neq 0$ then we can choose $w \in V^{\bot} \setminus \{0\}$ such that $$A \left \|w \right \|^{2} \leq \sum \limits_{l\in \mathfrak{L}} \sum \limits_{j \in \mathbb{Z}} \sum \limits_{a \in I_{p}}| \langle w,f^{(l)}_{j,a} \rangle |^{2} =0.$$ This imply $w=0$ which is a contradiction. Hence $L^{2}(\mathbb{Q}_{p})=V$. Thus a multiframelet in $p$-adic setting is complete.

\begin{remark}
	Using Proposition \ref{lower}, it has been verified that every multiframelet is an extension of orthonormal basis.
\end{remark}

\begin{definition}(Multiframelet Sequence). A sequence $\{f_{j,a}^{(l)} : j \in \mathbb{Z},~ a \in I_{p},~ l \in \mathfrak{L}\}$ is said to be a multiframelet sequence if it is a multiframelet for $\overline{\text{span}} \{f_{j,a}^{(l)} : j \in \mathbb{Z},~ a \in I_{p},~ l \in \mathfrak{L} \}$.
\end{definition}

\begin{proposition}
	Let $\{f_{j,a}^{(l)} : j \in \mathbb{Z},~ a \in I_{p},~ l \in \mathfrak{L} \}$ be a multiframelet sequence with pre-multiframelet operator $\mathcal{T}$. Then $\mathcal{R}_{\mathcal{T}^{*}\mathcal{T}} = \mathcal{R}_{\mathcal{T}^*}$.
\end{proposition}

\noindent \textbf{Proof :}
	Clearly, $\mathcal{T}^{*}\mathcal{T}$ maps $\mathcal{R}_{\mathcal{T}^*}$ into itself. $\{f_{j,a}^{(l)} : j \in \mathbb{Z},~ a \in I_{p},~ l \in \mathfrak{L}\}$ is a multiframelet sequence imply $\mathcal{R}_{\mathcal{T}} = \overline{span} \{f_{j,a}^{(l)} : j \in \mathbb{Z},~ a \in I_{p},~ l \in \mathfrak{L} \}$. As $L^{2}(\mathbb{Q}_{p}) = \mathcal{R}_{\mathcal{T}} \bigoplus \mathcal{R}_{\mathcal{T}}^{\perp} = R_{\mathcal{T}} \bigoplus \mathcal{N}_{\mathcal{T}^*}$, any $g \in L^{2}(\mathbb{Q}_{p})$ can be written as
	\begin{equation*}
	g = \mathcal{T} \{ c(j,a,l) : j \in \mathbb{Z},~ a \in I_{p},~ l \in \mathfrak{L} \} + h,
	\end{equation*}
	where $c(l,j,a) \in \ell^{2}(\mathfrak{L} \times \mathbb{Z} \times I_{p})$ and $h \in \mathcal{N}_{T^*}$. Then
	\begin{equation*}
	\mathcal{T}^{*}g = \mathcal{T}^{*}\mathcal{T} \{ c(l,j,a) ~:~ j \in \mathbb{Z},~ a \in I_{p},~ l \in \mathfrak{L} \}.
	\end{equation*}
	Therefore $\mathcal{R}_{\mathcal{T}^{*}\mathcal{T}} = \mathcal{R}_{\mathcal{T}^*}$.

Image of a multiframelet sequence under unitary operator is also a multiframelet sequence, which is evident from the following theorem. 

\begin{theorem}
	Let $\{f_{j,a}^{(l)} : j \in \mathbb{Z},~ a \in I_{p},~ l \in \mathfrak{L} \}$ be a multiframelet sequence with bounds $A,~B$ and $\mathcal{U} : L^{2}(\mathbb{Q}_{p}) \to L^{2}(\mathbb{Q}_{p})$ is an unitary operator. Then $\{\mathcal{U} f_{j,a}^{(l)} : j \in \mathbb{Z},~ a \in I_{p},~ l \in \mathfrak{L} \}$ is a multiframelet sequence with same bounds.
\end{theorem}

\noindent \textbf{Proof :}
	Let $g \in span \{\mathcal{U} f_{j,a}^{(l)} : j \in \mathbb{Z},~ a \in I_{p},~ l \in \mathfrak{L} \}$. Then there is $h \in span \{ f_{j,a}^{(l)} : j \in \mathbb{Z},~ a \in I_{p},~ l \in \mathfrak{L} \}$ such that $g=\mathcal{U}h$. As $\mathcal{U}$ is unitary, we have $$\sum \limits_{l\in \mathfrak{L}} \sum \limits_{j \in \mathbb{Z}} \sum \limits_{a \in I_{p}}| \langle g, \mathcal{U}f^{(l)}_{j,a} \rangle |^{2}
	= \sum \limits_{l\in \mathfrak{L}} \sum \limits_{j \in \mathbb{Z}} \sum \limits_{a \in I_{p}}| \langle \mathcal{U}h, \mathcal{U}f^{(l)}_{j,a} \rangle |^{2}
	= \sum \limits_{l\in \mathfrak{L}} \sum \limits_{j \in \mathbb{Z}} \sum \limits_{a \in I_{p}}| \langle h, f^{(l)}_{j,a} \rangle |^{2}.$$ Since $\{f_{j,a}^{(l)} : j \in \mathbb{Z},~ a \in I_{p},~ l \in \mathfrak{L} \}$ is a multiframelet sequence with bounds $A,~B; \\
	~~~~~ A \|h\|^{2} \leq \sum \limits_{l\in \mathfrak{L}} \sum \limits_{j \in \mathbb{Z}} \sum \limits_{a \in I_{p}}| \langle h, f^{(l)}_{j,a} \rangle |^{2} \leq B \|h\|^{2} \\
	~~~~~ \text{i.e.} ~~ A \| \mathcal{U}h \|^{2} \leq \sum \limits_{l\in \mathfrak{L}} \sum \limits_{j \in \mathbb{Z}} \sum \limits_{a \in I_{p}}| \langle \mathcal{U}^{-1}g, f^{(l)}_{j,a} \rangle |^{2} \leq B \| \mathcal{U}h \|^{2},~~~~ (\mathcal{U}~ \text{is unitary} \Rightarrow \|\mathcal{U}h\| = \|h\|)  \\
	~~~~~ \text{i.e.}~~ A \|g\|^{2} \leq \sum \limits_{l\in \mathfrak{L}} \sum \limits_{j \in \mathbb{Z}} \sum \limits_{a \in I_{p}}| \langle \mathcal{U}^{*}g, f^{(l)}_{j,a} \rangle |^{2} \leq B \|g\|^{2},~~~~ (\mathcal{U}~ \text{is unitary} \Rightarrow \mathcal{U}^{*} = \mathcal{U}^{-1})\\
	~~~~~ \text{i.e.}~~ A \|g\|^{2} \leq \sum \limits_{l\in \mathfrak{L}} \sum \limits_{j \in \mathbb{Z}} \sum \limits_{a \in I_{p}}| \langle g, \mathcal{U}f^{(l)}_{j,a} \rangle |^{2} \leq B \|g\|^{2}.$\\
	By Proposition 3.3 of \cite{HB19}, above inequality holds for $g \in \overline{span} \{\mathcal{U} f_{j,a}^{(l)} : j \in \mathbb{Z},~ a \in I_{p},~ l \in \mathfrak{L} \}.$ Therefore $\{\mathcal{U} f_{j,a}^{(l)} : j \in \mathbb{Z},~ a \in I_{p},~ l \in \mathfrak{L} \}$ is a multiframelet sequence with bounds $A,~B$.

\section{Multiframelet Operator and its Dual}

By composing synthesis and analysis operator, we obtain the multiframelet operator\\ $\mathcal{S} :L^{2}(\mathbb{Q}_{p}) \rightarrow L^{2}(\mathbb{Q}_{p})$ defined by $\mathcal{S}g=\mathcal{T}\mathcal{T}^{*}g= \sum \limits_{l\in \mathfrak{L}} \sum \limits_{j \in \mathbb{Z}} \sum \limits_{a \in I_{p}} \langle g,f^{(l)}_{j,a} \rangle f^{(l)}_{j,a}$. Now Equation \eqref{fie} can be rewritten as $A\left \| g \right|\|^{2} \leq \langle \mathcal{S}g, g \rangle \leq B\left \| g \right \|^{2}, ~~ \forall ~ g \in L^{2}(\mathbb{Q}_{p})$. 

Duffin and Schaeffer \cite{1} have studied properties of frame operator in $\mathbb{R}$. Heil \cite{11} continued to study this in general Hilbert space. Later, Debnath \cite{14} also independently studied this for general Hilbert space. Here we have shown similar results also hold in $p$-adic setting. We discuss some important properties of multiframelet operator. In this context, Haldar and Bhandari \cite{HB19} has proved next two important results of multiframelet operator and its inverse.

\begin{proposition}
	$\mathcal{S}$ is a linear, bijective, self-adjoint, positive, bounded operator.
\end{proposition}

\begin{proposition}
	If $\{f^{(l)}_{j,a} : j \in \mathbb{Z},~ a \in I_{p},~ l \in \mathfrak{L} \}$ is a multiframelet with bounds $A,~ B$, then $\{\mathcal{S}^{-1}f^{(l)}_{j,a} : j \in \mathbb{Z},~ a \in I_{p},~ l \in \mathfrak{L} \}$ is a multiframelet with bounds $B^{-1} ,~ A^{-1}$. Furthermore, if $A,~B$ are optimal for $\{ f^{(l)}_{j,a} : j \in \mathbb{Z},~ a \in I_{p},~ l \in \mathfrak{L} \}$ then $B^{-1} ,~ A^{-1}$ are optimal for $\{\mathcal{S}^{-1}f^{(l)}_{j,a} : j \in \mathbb{Z},~ a \in I_{p},~ l \in \mathfrak{L} \}$ whose multiframelet operator is $\mathcal{S}^{-1}$.
\end{proposition}

$\{\mathcal{S}^{-1}f^{(l)}_{j,a} : j \in \mathbb{Z},~ a \in I_{p},~ l \in \mathfrak{L} \}$ is called the canonical dual multiframelet of $\{f^{(l)}_{j,a} : j \in \mathbb{Z},~ a \in I_{p},~ l \in \mathfrak{L} \}$. Note that $\{\mathcal{S}f^{(l)}_{j,a} : j \in \mathbb{Z},~ a \in I_{p},~ l \in \mathfrak{L} \}$ is also a multiframelet.

\begin{theorem}
	Let $\mathcal{S}^{- \frac{1}{2}}$ denotes positive square root of $\mathcal{S}^{-1}$. Then $\{\mathcal{S}^{- \frac{1}{2}} f^{(l)}_{j,a} : j \in \mathbb{Z},~ a \in I_{p},~ l \in \mathfrak{L} \}$ is a normalized tight multiframe and for all $g \in L^{2}(\mathbb{Q}_{p})$,
	\begin{equation*}
	g = \sum \limits_{l\in \mathfrak{L}} \sum \limits_{j \in \mathbb{Z}} \sum \limits_{a \in I_{p}} \langle g ~,~ \mathcal{S}^{- \frac{1}{2}}f_{j,a}^{(l)} \rangle \mathcal{S}^{- \frac{1}{2}} f_{j,a}^{(l)}.
	\end{equation*}
\end{theorem}

\noindent \textbf{Proof :}
	As $\mathcal{S}^{- \frac{1}{2}}$ is defined as limit of a  sequence of polynomials in $\mathcal{S}^{-1}$, it commutes with $\mathcal{S}^{-1}$ and hence with $\mathcal{S}$ (cf. \cite{6}). Again we have,
	\begin{eqnarray}
	g = \mathcal{S}^{- \frac{1}{2}} \mathcal{S}\mathcal{S}^{- \frac{1}{2}} g \nonumber
	= \mathcal{S}^{- \frac{1}{2}} \sum \limits_{l\in \mathfrak{L}} \sum \limits_{j \in \mathbb{Z}} \sum \limits_{a \in I_{p}} \langle \mathcal{S}^{- \frac{1}{2}}g ~,~ f_{j,a}^{(l)} \rangle  f_{j,a}^{(l)}
	= \sum \limits_{l\in \mathfrak{L}} \sum \limits_{j \in \mathbb{Z}} \sum \limits_{a \in I_{p}} \langle g ~,~ \mathcal{S}^{- \frac{1}{2}} f_{j,a}^{(l)} \rangle  \mathcal{S}^{- \frac{1}{2}} f_{j,a}^{(l)}.
	\end{eqnarray}
	Now let us consider the  following inner product
	\begin{eqnarray}
	\|g\|^{2} = \langle g,g \rangle \nonumber
	&=& \left \langle g ~, ~ \sum \limits_{l\in \mathfrak{L}} \sum \limits_{j \in \mathbb{Z}} \sum \limits_{a \in I_{p}} \langle g ~,~ \mathcal{S}^{- \frac{1}{2}} f_{j,a}^{(l)} \rangle  \mathcal{S}^{- \frac{1}{2}} f_{j,a}^{(l)} ~ \right \rangle \nonumber \\
	&=& \sum \limits_{l\in \mathfrak{L}} \sum \limits_{j \in \mathbb{Z}} \sum \limits_{a \in I_{p}} \overline{ \left \langle g ~,~ \mathcal{S}^{- \frac{1}{2}} f_{j,a}^{(l)} \right \rangle} \left \langle g ~, ~ \mathcal{S}^{- \frac{1}{2}} f_{j,a}^{(l)} ~ \right \rangle \nonumber \\
	&=& \sum \limits_{l\in \mathfrak{L}} \sum \limits_{j \in \mathbb{Z}} \sum \limits_{a \in I_{p}}  \left |\left \langle g ~, ~ \mathcal{S}^{- \frac{1}{2}} f_{j,a}^{(l)} ~ \right \rangle \right |^{2}. \nonumber
	\end{eqnarray}
	This completes the proof.

\begin{proposition}
	Let $G$ be the analysis operator associated with the multiframelet $\{f^{(l)}_{j,a} : j \in \mathbb{Z},~ a \in I_{p},~ l \in \mathfrak{L} \}$ and $\widetilde{G}$ be the analysis operator corresponding to its canonical dual multiframelet $\{S^{-1}f^{(l)}_{j,a} : j \in \mathbb{Z},~ a \in I_{p},~ l \in \mathfrak{L} \}$. Then $\widetilde{G}^{*} G = \mathcal{I} = G^{*} \widetilde{G}$.
\end{proposition}

\noindent \textbf{Proof :} We have
	\begin{eqnarray}
	G(G^{*}G)^{-1}g
	&=& \{ \langle (G^{*}G)^{-1}g ~,~ f_{j,a}^{(l)} \rangle : j \in \mathbb{Z},~ a \in I_{p},~ l \in \mathfrak{L} \} \nonumber \\
	&=& \{ \langle g ~,~ (G^{*}G)^{-1}f_{j,a}^{(l)} \rangle : j \in \mathbb{Z},~ a \in I_{p},~ l \in \mathfrak{L} \} \nonumber \\
	&=& \{ \langle g ~,~ S^{-1}f_{j,a}^{(l)} \rangle : j \in \mathbb{Z},~ a \in I_{p},~ l \in \mathfrak{L} \} \nonumber \\
	&=& \widetilde{G}g \nonumber
	\end{eqnarray}
	So, $G(G^{*}G)^{-1}= \widetilde{G}$. Now $\widetilde{G}^{*}G = (G(G^{*}G)^{-1})^{*}G = (G^{*}G)^{-1} G^{*}G = \mathcal{I}$ and also $G^{*} \widetilde{G} = G^{*}G(G^{*}G)^{-1} = \mathcal{I}$. Hence $\widetilde{G}^{*} G = \mathcal{I} = G^{*} \widetilde{G}$.

The notion of multiframelet presents decomposition technique, every element in $L^{2}(\mathbb{Q}_{p})$ can be written as infinite linear combination of its multiframelets. So we can think multiframelet as a kind of basis, though it is already proved by Haldar and Bhandari \cite{HB19}, here we presents an alternative proof of same.  

\begin{theorem} \label{md}(Multiframelet decomposition).
	Let $\mathfrak{f}$=$\{f^{(1)}, \ldots, f^{(L)}\}$ is a multiframelet for $L^{2}(\mathbb{Q}_{p})$ with multiframelet operator $\mathcal{S}$. Then any $g \in L^{2}(\mathbb{Q}_{p})$ can be written as
	\begin{equation*}
	g = \sum \limits_{l\in \mathfrak{L}} \sum \limits_{j \in \mathbb{Z}} \sum \limits_{a \in I_{p}} \langle g ,~ \mathcal{S}^{-1} f_{j,a}^{(l)} \rangle f_{j,a}^{(l)} = \sum \limits_{l\in \mathfrak{L}} \sum \limits_{j \in \mathbb{Z}} \sum \limits_{a \in I_{p}} \langle g ,~ f_{j,a}^{(l)} \rangle \mathcal{S}^{-1} f_{j,a}^{(l)},
	\end{equation*}
	where sum converges unconditionally.
\end{theorem}

\noindent \textbf{Proof :}
	Let $G$ be the analysis operator associated to multiframelet $\{f^{(l)}_{j,a} : j \in \mathbb{Z},~ a \in I_{p},~ l \in \mathfrak{L} \}$ and let $\widetilde{G}$ be the analysis operator corresponding to its dual multiframelet $\{\mathcal{S}^{-1}f^{(l)}_{j,a} : j \in \mathbb{Z},~ a \in I_{p},~ l \in \mathfrak{L} \}$. Since $\mathcal{I} = G^{*} \widetilde{G}, ~\forall ~ g \in L^{2}(\mathbb{Q}_{p})$, we have $g = G^{*} \widetilde{G} g = G^{*} \{ \langle g ,~ \mathcal{S}^{-1}f_{j,a}^{(l)} \rangle : j \in \mathbb{Z},~ a \in I_{p},~ l \in \mathfrak{L} \} = \sum \limits_{l\in \mathfrak{L}} \sum \limits_{j \in \mathbb{Z}} \sum \limits_{a \in I_{p}} \langle g ,~ \mathcal{S}^{-1} f_{j,a}^{(l)} \rangle f_{j,a}^{(l)}$.\\
	The proof of the other equality is similar and needs to use $\mathcal{I} = \widetilde{G}^{*} G$. 

\begin{remark}
	The numbers $\langle g ,~ \mathcal{S}^{-1} f_{j,a}^{(l)} \rangle$ are called multiframelet coefficients. These are very important as these contains information of decomposed function $g$. Theorem \ref{md} also shows that any element in $L^{2}(\mathbb{Q}_{p})$ can be written in-terms of canonical dual multiframelet. Furthermore, $\mathcal{S}$ is surjective and therefore a topological isomorphism of $L^{2}(\mathbb{Q}_{p})$. When $\mathfrak{f}$ is a tight multiframelet then $\mathcal{S}^{-1}=A^{-1}~ \mathcal{I}$ and representation of $g$ in Theorem \ref{md} changes to $g = A^{-1} \sum \limits_{l\in \mathfrak{L}} \sum \limits_{j \in \mathbb{Z}} \sum \limits_{a \in I_{p}} \langle g ,~ f_{j,a}^{(l)} \rangle f_{j,a}^{(l)}$. It is to be noted that $A$ is an eigen value of $\mathcal{S}$ in that case. Also natural question one may ask, Is one multiframelet has representation by other multiframelet in $L^{2}(\mathbb{Q}_{p})$ ? The answer is affirmative by multiframelet decomposition theorem.
\end{remark}

The following proposition shows that scalars given in Theorem \ref{md} has minimum $\ell^{2}$-norm among all choices of scalars in $g= \sum \limits_{l\in \mathfrak{L}} \sum \limits_{j \in \mathbb{Z}} \sum \limits_{a \in I_{p}} d(l,j,a) f_{j,a}^{(l)}$.

\begin{proposition}
	Let $\mathfrak{f}$=$\{f^{(1)}, \ldots, f^{(L)}\}$ be a multiframelet and $g \in L^{2}(\mathbb{Q}_{p})$. If $g= \sum \limits_{l=1}^{L} \sum \limits_{j \in \mathbb{Z}} \sum \limits_{a \in I_{p}} d(l,j,a) f_{j,a}^{(l)}$ then following equality holds
	\begin{equation*}
	\sum \limits_{l\in \mathfrak{L}} \sum \limits_{j \in \mathbb{Z}} \sum \limits_{a \in I_{p}} |d(l,j,a)|^{2}
	= \sum \limits_{l\in \mathfrak{L}} \sum \limits_{j \in \mathbb{Z}} \sum \limits_{a \in I_{p}} |\langle g,~ \mathcal{S}^{-1}f_{j,a}^{(l)} \rangle |^{2} + \sum \limits_{l\in \mathfrak{L}} \sum \limits_{j \in \mathbb{Z}} \sum \limits_{a \in I_{p}}  |d(l,j,a) - \langle g,~ \mathcal{S}^{-1}f_{j,a}^{(l)} \rangle |^{2}.
	\end{equation*}
\end{proposition}

\noindent \textbf{Proof :} By Theorem \ref{md}, we have $g = \sum \limits_{l\in \mathfrak{L}} \sum \limits_{j \in \mathbb{Z}} \sum \limits_{a \in I_{p}} \langle g ,~ \mathcal{S}^{-1} f_{j,a}^{(l)} \rangle f_{j,a}^{(l)}$. Again for every $g \in L^{2}(\mathbb{Q}_{p})$ we have,
	\begin{eqnarray}
	\langle g,~ \mathcal{S}^{-1} g \rangle
	&=& \left \langle \sum \limits_{l\in \mathfrak{L}} \sum \limits_{j \in \mathbb{Z}} \sum \limits_{a \in I_{p}} d(l,j,a) f_{j,a}^{(l)} ~,~ \mathcal{S}^{-1} g \right \rangle  \nonumber\\
	&=&   \sum \limits_{l\in \mathfrak{L}} \sum \limits_{j \in \mathbb{Z}} \sum \limits_{a \in I_{p}} d(l,j,a) \left \langle f_{j,a}^{(l)} ~,~ \mathcal{S}^{-1} g \right \rangle  \nonumber\\
	&=&   \sum \limits_{l\in \mathfrak{L}} \sum \limits_{j \in \mathbb{Z}} \sum \limits_{a \in I_{p}} d(l,j,a) \left \langle \mathcal{S}^{-1} f_{j,a}^{(l)} ~,~ g \right \rangle  \nonumber\\
	&=& \sum \limits_{l\in \mathfrak{L}} \sum \limits_{j \in \mathbb{Z}} \sum \limits_{a \in I_{p}} d(l,j,a) ~ \overline{\left \langle g,~ \mathcal{S}^{-1} f_{j,a}^{(l)} \right \rangle} \nonumber
	\end{eqnarray}
	
	\begin{eqnarray}
	\text{also,  }~~~ \langle g,~ \mathcal{S}^{-1} g \rangle
	&=& \left \langle \sum \limits_{l\in \mathfrak{L}} \sum \limits_{j \in \mathbb{Z}} \sum \limits_{a \in I_{p}} \langle g ,~ \mathcal{S}^{-1} f_{j,a}^{(l)} \rangle f_{j,a}^{(l)} ~,~ \mathcal{S}^{-1} g \right \rangle  \nonumber\\
	&=&  \sum \limits_{l\in \mathfrak{L}} \sum \limits_{j \in \mathbb{Z}} \sum \limits_{a \in I_{p}} \langle g ,~ \mathcal{S}^{-1} f_{j,a}^{(l)} \rangle ~ \langle f_{j,a}^{(l)} ~,~ \mathcal{S}^{-1} g \rangle  \nonumber\\
	&=&  \sum \limits_{l\in \mathfrak{L}} \sum \limits_{j \in \mathbb{Z}} \sum \limits_{a \in I_{p}} \langle g ,~ \mathcal{S}^{-1} f_{j,a}^{(l)} \rangle ~ \langle \mathcal{S}^{-1}f_{j,a}^{(l)} ~,~ g  \rangle  \nonumber\\
	&=&  \sum \limits_{l\in \mathfrak{L}} \sum \limits_{j \in \mathbb{Z}} \sum \limits_{a \in I_{p}} \langle g ,~ \mathcal{S}^{-1} f_{j,a}^{(l)} \rangle ~ \overline{\langle g,~ \mathcal{S}^{-1}f_{j,a}^{(l)} \rangle} \nonumber
	\end{eqnarray}
	So $\{d(l,j,a) - \langle g,~ \mathcal{S}^{-1}f_{j,a}^{(l)} \rangle : j \in \mathbb{Z},~ a \in I_{p},~ l \in \mathfrak{L} \}$ is orthogonal with $\{\langle g,~ \mathcal{S}^{-1}f_{j,a}^{(l)} \rangle : j \in \mathbb{Z},~ a \in I_{p},~ l \in \mathfrak{L} \}$ in $\ell^{2}$-norm. Therefore
	$$\| \{d(l,j,a)\}\|^{2} = \|\{\langle g,~ \mathcal{S}^{-1}f_{j,a}^{(l)} \rangle \}\|^{2} + \|\{d(l,j,a) - \langle g,~ \mathcal{S}^{-1}f_{j,a}^{(l)} \rangle \}\|^{2}$$
	
	$$\text{i.e. } \sum \limits_{l\in \mathfrak{L}} \sum \limits_{j \in \mathbb{Z}} \sum \limits_{a \in I_{p}} |d(l,j,a)|^{2} = \sum \limits_{l\in \mathfrak{L}} \sum \limits_{j \in \mathbb{Z}} \sum \limits_{a \in I_{p}} |\langle g,~ \mathcal{S}^{-1}f_{j,a}^{(l)} \rangle |^{2} + \sum \limits_{l\in \mathfrak{L}} \sum \limits_{j \in \mathbb{Z}} \sum \limits_{a \in I_{p}}  |d(l,j,a) - \langle g,~ \mathcal{S}^{-1}f_{j,a}^{(l)} \rangle |^{2}.$$

\begin{proposition}
	Let $\{f_{j,~a}^{(l)} : j \in \mathbb{Z},~ a \in I_{p},~ l \in \mathfrak{L} \}$ be a multiframelet sequence in $L^{2}(\mathbb{Q}_{p})$. Then the orthogonal projection $\mathfrak{D}$ of  $\ell^{2}(\mathfrak{L} \times \mathbb{Z} \times I_{p})$ onto $\mathcal{R}_{\mathcal{T}^*}$ is given by
	\begin{eqnarray*}
		\mathfrak{D} \left \{c(l,j,a) : l \in \mathfrak{L}, j \in \mathbb{Z}, a \in I_{p}\right \} = \left \{ \left \langle \sum \limits_{l\in \mathfrak{L}} \sum \limits_{j \in \mathbb{Z}} \sum \limits_{a \in I_{p}} c(l,j,a)~ \mathcal{S}^{-1}f_{j,a}^{(l)}~,~ f_{j^{\prime},a^{\prime}}^{(l^{\prime})} \right \rangle : j^{\prime} \in \mathbb{Z},~ a^{\prime} \in I_{p},~ l^{\prime} \in \mathfrak{L} \right \}.
	\end{eqnarray*}
\end{proposition}

\noindent \textbf{Proof :}
	It is enough to show that $\mathfrak{D}$ is identity map on $\mathcal{R}_{\mathcal{T}^*}$ and zero map on $\mathcal{R}_{\mathcal{T}^*}^{\perp}$.
	\begin{eqnarray*} 
		\text{Now, } &&\mathfrak{D} \left \{ \left \langle g,~f_{j,a}^{(l)} \right \rangle  :~ l \in \mathfrak{L}, j \in \mathbb{Z}, a \in I_{p} \right \} \\
		&&= \left \{ \left \langle \sum \limits_{l\in \mathfrak{L}} \sum \limits_{j \in \mathbb{Z}} \sum \limits_{a \in I_{p}} \langle g,~f_{j,a}^{(l)} \rangle \mathcal{S}^{-1} f_{j,a}^{(l)}~,~ f_{j^{\prime},a^{\prime}}^{(l^{\prime})} \right \rangle : j^{\prime} \in \mathbb{Z},~ a^{\prime} \in I_{p},~ l^{\prime} \in \mathfrak{L} \right \} \\
		&&= \left \{ \left \langle g ~,~ f_{j^{\prime},a^{\prime}}^{(l^{\prime})} \right \rangle : j^{\prime} \in \mathbb{Z},~ a^{\prime} \in I_{p},~ l^{\prime} \in \mathfrak{L} \right \}
	\end{eqnarray*}
	
	So, $\mathfrak{D}$ is identity on $\mathcal{R}_{\mathcal{T}^*}$.
	As $\mathcal{R}_{\mathcal{T}^*}^{\perp} = \mathcal{N}_{\mathcal{T}}$, $$\{c(l,j,a) : l \in \mathfrak{L}, j \in \mathbb{Z}, a \in I_{p} \} \in \mathcal{N}_{\mathcal{T}} \Rightarrow \sum \limits_{l=1}^{L} \sum \limits_{j \in \mathbb{Z}} \sum \limits_{a \in I_{p}} c(l,j,a) f_{j,a}^{(l)} = 0.$$ 
	
	\begin{eqnarray*}
		\text{Now,} &&\mathfrak{D}\{c(l,j,a) : l \in \mathfrak{L}, j \in \mathbb{Z}, a \in I_{p}\} \\
		&&= \left \{ \left \langle \sum \limits_{l\in \mathfrak{L}} \sum \limits_{j \in \mathbb{Z}} \sum \limits_{a \in I_{p}} c(l,j,a) S^{-1} f_{j,a}^{(l)}~,~ f_{j^{\prime},a^{\prime}}^{(l^{\prime})} \right \rangle : j^{\prime} \in \mathbb{Z},~ a^{\prime} \in I_{p},~ l^{\prime} \in \mathfrak{L} \right \} \\
		&&= \left \{ \left \langle \mathcal{S}^{-1} \sum \limits_{l\in \mathfrak{L}} \sum \limits_{j \in \mathbb{Z}} \sum \limits_{a \in I_{p}} c(l,j,a)  f_{j,a}^{(l)}~,~ f_{j^{\prime},a^{\prime}}^{(l^{\prime})} \right \rangle : j^{\prime} \in \mathbb{Z},~ a^{\prime} \in I_{p},~ l^{\prime} \in \mathfrak{L} \right \} \\
		&&= \left \{ \left \langle \mathcal{S}^{-1} 0 ~,~ f_{j^{\prime},a^{\prime}}^{(l^{\prime})} \right \rangle : j^{\prime} \in \mathbb{Z},~ a^{\prime} \in I_{p},~ l^{\prime} \in \mathfrak{L} \right \}\\
		&&=  \{ 0 : j^{\prime} \in \mathbb{Z},~ a^{\prime} \in I_{p},~ l^{\prime} \in \mathfrak{L} \}
	\end{eqnarray*}
	
	Therefore $\mathfrak{D}$ is zero on $\mathcal{R}_{\mathcal{T}^*}^{\perp}$ and this completes the proof.

\section{Multiframelet Set}

In this section, we introduce multiframelet set and raise one important question about it.

\begin{definition} (Multiframelet set).
	A measurable set $\mathfrak{F} \subset \mathbb{Q}_{p}$ is said to be a multiframelet set of order $L$ if $\mathfrak{F}=\bigcup \limits_{l=1}^{L} F_{l}$ where $\hat{f}^{(l)}=\mathds{1}_{F_{l}}, ~~ \forall ~l$ and $\mathfrak{f}$=$\{f^{(1)}, \ldots, f^{(L)}\}$ is a multiframelet for $L^{2}(\mathbb{Q}_{p})$.
\end{definition}

\begin{example} $\mathfrak{F}_{p}=B_{0}(-\frac{1}{p}) \cup \ldots \cup B_{0}(-\frac{p-1}{p})$ ia a multiframelet set of order $(p-1)$.
\end{example}

\begin{remark}
	It is easy to see that $\mathfrak{F}_{p}$ is open, compact but not connected. We are yet to answer these properties for general multiframelet sets in $\mathbb{Q}_{p}$.
\end{remark}

\begin{theorem}
	If $\mathfrak{f}$=$\{f^{(1)}, \ldots, f^{(L)}\}$ is a normalized tight multiframelet for $L^{2}(\mathbb{Q}_{p})$ where $\hat{f^{(l)}}= \mathds{1}_{F_{l}} \text{ and } F_{l} \subset \mathbb{Q}_{p}, ~~ \forall ~l$. Then $$\|g\|^{2} = \sum \limits_{l\in \mathfrak{L}}  \sum \limits_{a \in I_{p}}  \sum \limits_{j \in \mathbb{Z}} p^{-j} \left |~ \int \limits_{p^{-j} F_{l}}  \chi_{p}(p^{j}a \xi)~ \overline{\hat{g}(\xi)} ~ d \xi \right |^{2}, ~~~ \forall g \in L^{2}(\mathbb{Q}_{p}).$$
\end{theorem}

\noindent \textbf{Proof :} As $\mathfrak{f}$ is normalized tight multiframelet, $\|g\|^{2} = \sum \limits_{l\in \mathfrak{L}} \sum \limits_{a \in I_{p}}  \sum \limits_{j \in \mathbb{Z}}~| \langle f_{j,a}^{(l)} ~, g  \rangle |^{2}, ~~~ \forall ~ g \in L^{2}(\mathbb{Q}_{p}).$
	
	Now using Perseval theorem, we get
	\begin{eqnarray*}
		\|g\|^{2} &=& \sum \limits_{l\in \mathfrak{L}} \sum \limits_{a \in I_{p}}  \sum \limits_{j \in \mathbb{Z}}~ | \langle \hat{f_{j,a}^{(l)}} ~, \hat{g} \rangle |^{2} \\
		&=& \sum \limits_{l\in \mathfrak{L}} \sum \limits_{a \in I_{p}}  \sum \limits_{j \in \mathbb{Z}}~ \left|~ \int \limits_{\mathbb{Q}_{p}} p^{-\frac{j}{2}}~ \chi_{p}(p^{j}a \xi) ~ \hat{f^{(l)}}(p^{j}\xi) ~ \overline{\hat{g}(\xi)} ~ d \xi \right|^{2} \\ &=& \sum \limits_{l\in \mathfrak{L}} \sum \limits_{a \in I_{p}}  \sum \limits_{j \in \mathbb{Z}}~ \left|~ \int \limits_{\mathbb{Q}_{p} ~ \cap ~ p^{-j}F_{l}} p^{-\frac{j}{2}} ~ \chi_{p}(p^{j}a \xi) ~ \overline{\hat{g}(\xi)} ~ d \xi \right|^{2} \\ &=& \sum \limits_{l\in \mathfrak{L}} \sum \limits_{a \in I_{p}}  \sum \limits_{j \in \mathbb{Z}}~ p^{-j} \left|~ \int \limits_{p^{-j}F_{l}}  ~ \chi_{p}(p^{j}a \xi) ~ \overline{\hat{g}(\xi)} ~ d \xi \right|^{2}
	\end{eqnarray*}
	Hence our assertion is tenable.

\section*{Acknowledgements}

This research was supported by MHRD, Government of India. Moreover, The author acknowledges insightful comments by  Animesh Bhandari to improve this article. He also thanks to Dr. Divya Singh for her help and encouragement.

\end{document}